\documentclass[12pt]{article}
\usepackage{amsmath,amssymb,amsthm,fullpage}
\usepackage[numbers]{natbib}

\newcommand{\R}{{\mathbb R}}

\newcommand{\ep}{{\varepsilon}}

\newcommand{\dist}{\text{dist}}

\usepackage{braket}

\newcommand{\parens}[1]{\left( #1 \right)}

\theoremstyle{plain}
\newtheorem{theorem}{Theorem}[section]

\newtheorem{lemma}[theorem]{Lemma}

\newtheorem{example}[theorem]{Example}

\theoremstyle{remark}
\newtheorem*{remark}{Remark}

\theoremstyle{definition}
\newtheorem{definition}{Definition}[section]

\usepackage{amsmath}
\usepackage{amssymb}
\usepackage{latexsym}
\usepackage{amsfonts,setspace}
\usepackage{fullpage}
\usepackage{amsthm}
\usepackage{tikz}
\usepackage{graphicx}
\graphicspath{ {./images/} }

\usepackage{color}
\definecolor{red}{rgb}{.8,0,0}
\definecolor{green}{rgb}{0,.7,0}
\definecolor{blue}{rgb}{0,0,.8}
\newcommand{\peter}[1]{\textbf{{\color{red}[Peter: #1]}}}
\title{On the Bipartiteness Constant and Expansion of Cayley Graphs}
\author{Nina Moorman \thanks{School of Mathematics and School of Computer Science, Georgia Institute of Technology.  Supported in part by the NSF grants DMS-1811935 and
NSF TRIPODS-1740776.} \and Peter Ralli \thanks{Department of Mathematics, Weizmann Institute of Science.  Supported in part by European Research Council grant no. 803048.  Part of the research was conducted while this author was a scientific collaborator at the \'{E}cole polytechnique f\'{e}d\'{e}rale de Lausanne.}\and Prasad Tetali \footnotemark[1] }
\date{\today}
\begin{document}

\maketitle

\begin{abstract}
    Let $G$ be a finite, undirected, $d$-regular graph and $A(G)$ its normalized adjacency matrix, with eigenvalues  $1 = \lambda_1(A)\geq \dots \ge \lambda_n \ge -1$. It is a classical fact that $\lambda_n = -1$ if and only if $G$ is bipartite. Our main result provides a quantitative separation of $\lambda_n$ from $-1$ in the case of Cayley graphs, in terms of their expansion.  Denoting $h_{out}$ by the (outer boundary) vertex expansion of $G$, we show that if $G$ is a non-bipartite Cayley graph (constructed using a group and a symmetric generating set of size $d$) then $\lambda_n \ge -1 + \frac{ch_{out}^2}{d^2}\,,$ for $c$ an absolute constant.  We exhibit graphs for which this result is tight up to a factor depending on $d$. This improves upon a recent result by Biswas and Saha \cite{BS} who showed $\lambda_n \ge -1 + \frac{h_{out}^4}{2^9d^8}\,.$ We also note that such a result could not be true for general  non-bipartite graphs. 
    
\end{abstract}
\section{Introduction}

It is well-known that the leading eigenvalue $\lambda_1$ of the normalized adjacency matrix of a regular graph is $1$.  A topic of much interest is the second-largest absolute eigenvalue, $\max(\lambda_2, -\lambda_n$).  The famous Cheeger inequalities relate the \textit{spectral gap} $1-\lambda_2$ to the isoperimetric constant (see e.g.,\cite{Alon, Chung}).  We will be more concerned with what is in some sense the ``other spectral gap", namely the gap between $\lambda_n$ and $-1$, as studied in \cite{BJ,Trevisan}.

In a recent work, Breuillard et al.~\cite{BGGT} argued that if a non-bipartite Cayley graph is a  combinatorial expander - in the sense that $h_{out}$ is bounded away from $0$ - then it must also be a  spectral ``expander" in the sense that $\lambda_n$ is bounded away from $-1$ (\cite{BGGT}, Proposition E.1).  Combining this result with the Cheeger inequality, it is seen that $\max_{i>1}|\lambda_i|$ is bounded away from $1$.  Biswas~(\cite{Biswas}, Theorem 1.4), building on that argument, gave a bound of the form $1+\lambda_n \geq \frac{h_{out}^4}{2^9(d+1)^2d^6}$, and in a very recent work Biswas and Saha~(\cite{BS}, Theorem 2.12) refined the bound to $1+\lambda_n \geq \frac{h_{out}^4}{2^9d^8}$.  In this paper we improve on these results (see  Theorem~\ref{thm_finalbound} below), by proving that for every non-bipartite Cayley graph, $$1+\lambda_n \geq \frac{Ch_{out}^2}{d^2},$$ where $C>0$ is a universal constant.

Trevisan~\cite{Trevisan}, and independently Bauer and Jost~\cite{BJ}, introduced $\beta$ (defined below), a combinatorial parameter which measures the fraction of edges contributing to the non-bipartiteness of a graph.  $\beta$ is also a modification of a related constant developed much earlier by Desai and Rao~\cite{DR}.  An analogue of the Cheeger inequalities relates $\beta$ to $\lambda_n$:
\begin{theorem}[Trevisan \cite{Trevisan}, Equation~(8)]\label{thm_trevisan}
For any regular graph, $$2\beta \geq (1+\lambda_n)\geq \tfrac{1}{2}\beta^2.$$
\end{theorem}

Recall that $\lambda_n=-1$ if and only if the graph is bipartite, and $\beta\ge 0$ which was introduced by \cite{Trevisan} to capture non-bipartiteness, is also zero whenever the graph is bipartite. 

One can think of $\beta$ as serving the same role regarding $\lambda_n$ as the isoperimetric constant $h$ does for $\lambda_2$.  Following this analogy, we will define the \textit{outer vertex bipartiteness constant} $\beta_{out}$, just as $h_{out}$ is the outer (vertex) boundary isoperimetric constant.  In Theorem \ref{thm_betabounds} we demonstrate simple bounds relating $\beta_{out}$ to $\beta$, and by extension to $\lambda_n$.

\medskip

A brief outline of the overall proof strategies is as follows. The proofs of the previous results \cite{Biswas,BGGT} involve, for the Cayley graph $G = (X,S)$, examining the multigraph $G^2 = (X,S^2)$ with edges consisting of $2$-walks in $G$.  Then letting $A$ be the $h_{out}(G^2)$-achieving set, the outer (vertex) boundary $\partial_{out,G^2}(A)$ of $A$ in $G^2$ is bounded above by a function of $\lambda_n(G)$.  Following a method introduced by Freiman \cite{Freiman}, it is observed that if $\partial_{out,G^2} A$ is sufficiently small (with respect to $h_{out}(G))$, then there is a bipartition of $G$ which approximates $\{A,A^C\}$, contradicting the assumption that $G$ is not bipartite.

\medskip

Our innovation to this method is, rather than an $A$ as above, to consider sets $L,R$ which achieve $\beta_{out}(G)$; i.e., the best almost-bipartition of $G$.  In the proof of Theorem \ref{thm_vertexbound}, we will demonstrate an upper bound for $\partial_{out} L$ (or $\partial_{out} R$) in terms of $\beta_{out}$, and then, following the same method of Freiman, we argue that if $\beta_{out}$ is sufficiently small with respect to $h_{out}$ then $\{L,R\}$ approximates an actual bipartition of $G$, which gives a contradiction.  Our main result in Theorem \ref{thm_finalbound} follows by combining this proof with Trevisan's above-mentioned lower bound on $\lambda_n$.

\medskip

We demonstrate as Example \ref{example_cycle} that for the odd cycle $1+\lambda_n = \Theta(h_{out}^2)$, our result is tight up to a factor depending only on $d$, and therefore the term $h_{out}^2$ in our main result cannot be improved.
With Example \ref{example_no_converse} we demonstrate that there is no converse to our theorem; that is, there is no lower bound for $h_{out}$ in terms of $\lambda$ (or $\beta$) and $d$ even in the special case of non-bipartite Cayley graphs.

\subsection{Notations} Recall that a $d$-regular graph is one in which each vertex has exactly $d$ neighbors.  Throughout this paper we will let $G=(V,E)$ be a $d$-regular, connected, non-bipartite simple graph on $n$ vertices.
For a subset $S \subset V$, let $\partial(S)$ denote the edge boundary of $S$, namely the set of edges with precisely one endpoint in $S$.  $\partial_{out}(S)$ is the outer vertex  boundary, the set of vertices that are not in $S$ but do have a neighbor in $S$.
Let $h(G)$ denote the (edge) Cheeger constant of $G$ defined as 
$$h(G) = \min_{S\subset V\atop {0<|S|\le |V|/2}}\frac{|\partial(S)|}{d|S|}\,.$$  The classical isoperimetric constant {\em expansion} is defined using the outer vertex boundary: $$h_{out}(G) = \min_{S\subset V\atop {0<|S|\le |V|/2}}\frac{|\partial_{out}(S)|}{|S|}\,.$$

Let $A = [\tfrac{1}{d}1_{x\sim y}]$ be the normalized adjacency matrix and $\Delta = I-A$ the normalized Laplacian.  Write $1 = \lambda_1(A)\geq \dots \ge \lambda_n =: -1+\mu$.  Observe that $\mu >0$, since $G$ is assumed to be non-bipartite.

For a set $S\subset V$ we denote by $S^C$ its complement $V-S$.  For a pair of disjoint sets $L,R\subset V$,   define the bipartiteness ratio of $L$ and $R$ to be $$b(L,R) = \frac{e(L,R^C)+e(R,L^C)}{d|L\cup R|} = \frac{e(L,L)+e(R,R)+|\partial (L\cup R)|}{d|L\cup R|}.$$
Note: Our convention is that $e(A,B)$ counts the ordered pairs in $E(A,B):=\{(a,b)\in A\times B:a\sim b\}$, so that $e(A,A)$ equals twice the number of edges in the subgraph induced by $A$.

Following Trevisan \cite{Trevisan}, the {\em bipartiteness constant} of $G$ is $\beta(G) = \min_{L,R} b(L,R)$. As mentioned above, for a bipartite graph $G$, $\beta(G)=0$, since $L$ and $R$ can be chosen to be the bipartition of the graph, giving the numerator in $b(L,R)$ to be zero.
\medskip

One may also observe that the  definition of $\beta$ can be obtained by restricting the minimization in (the variational definition of) $\lambda_n$ (or rather that of $1+\lambda_n$) $$1+\lambda_n = \min_{x\in \R^n:x\neq \vec{0}}\frac{\tfrac{1}{d}\sum_{\{i,j\}\in E}|x_i+x_j|}{\sum_{i}|x_i|}$$ to functions taking values in $\{-1,0,+1\}$.

\medskip

In this work, we define the Cayley graph $(X,S)$, with $X$ being a finite group and $S$ a generating set of $X$, to have edges $g\sim gs$ for all $g\in X, s\in S$.


\section{Results}

In this section we investigate the problem of relating the bipartiteness constant $\beta$ to the isoperimetric constant $h$.  First, we will give a simple example that shows there cannot be a relationship in the general case of non-bipartite regular graphs.  But we find that there is no such obstruction for non-bipartite Cayley graphs, and those will be our main focus.

\subsection{General results}

There is no universal lower bound for $\beta$ (and, because $\mu \geq \tfrac{1}{2}\beta^2$, neither is there a bound for $\mu$) in terms of $h$ and $d$.  To see this we give an example that cannot follow any such bound.
\begin{example}
There is a non-bipartite graph for which $h$ is constant and $\beta = O(\tfrac{1}{nd})$.
\end{example}
\begin{proof}
Suppose $G$ is a bipartite expander, so that $h = c\in (0,1)$, $\beta = 0$.  Let $L,R$ be the bipartition of $G$ and let $(l_1,r_1),(l_2,r_2)\in L\times R$ be edges of $G$.  Now create the $d$-regular graph $G^*$ by replacing these two edges with $(l_1,l_2)$ and $(r_1,r_2)$.  Then $\beta(G^*) = \tfrac{4}{nd}$, achieved by the sets $L,R$.   Because we only changed a constant set of edges, $h(G^*) = h(G) + o_d(1) = c +o_d(1)$; i.e., $\beta$ is the smallest possible value and $h$ is constant.
\end{proof}

We will try to apply a similar idea to find an obstruction for Cayley graphs.

\begin{example}
There is a non-bipartite Cayley graph for which $h$ is constant and $\beta = O(\tfrac{1}{d})$.
\end{example}
\begin{proof}
Let $G$ be a bipartite Cayley graph for group $X$ with self-invertible set $S$, where $|S| = d$, and assume $h(G) = c\in (0,1)$.  Let $g\in X$ be an element in the same side of the bipartition as $e$.  Define $G^*$ to be the Cayley graph on $X$ generated by $S\cup \{g,g^{-1}\}$.  Then $\beta(G^*) =  \frac{2}{d+2}$ (or $\beta(G^*) = \frac{1}{d+1}$ in the case $g = g^{-1}$), which is achieved by $L,R$, while $h(G^*) \leq \frac{cd}{d+2}$.
\end{proof}

This example does not give as strong an obstruction, because we might have $\beta(G^*) = \Theta(\tfrac{1}{d})$, then $\beta(G^*)= \Omega(\frac{h(G^*)}{d})$; however a theorem stating that $\beta = \Omega(\frac{h}{d})$ for general Cayley graphs would certainly be of interest.  In the next section we will investigate what bounds are possible for this problem.

\subsection{Results for Cayley graphs}

For the proofs in this section it will be useful to define a bipartiteness parameter $\beta_{out}$ which involves a count of vertices that violate bipartiteness.  This is in contrast to the definition of $\beta$, which uses counts of edges.

\begin{definition}
For disjoint sets $L,R$, we define $$b_{out}(L,R) = \frac{I(L)+I(R)+|\partial_{out}(L\cup R)|}{|L\cup R|},$$ where $I(S)$ is the number of vertices in $S$ with a neighbor also in $S$.  The \em{outer vertex bipartiteness constant} is $$\beta_{out}(G) = \min_{L,R}b_{out}(L,R).$$
\end{definition}

As with $h$ and $h_{out}$, there is a simple relationship between $\beta$ and $\beta_{out}$.

\begin{theorem} \label{thm_betabounds}
$\beta_{out}(G)\geq \beta(G)\geq \tfrac{1}{d} \beta_{out}(G)$.
\end{theorem}
\begin{proof}
Take $L,R$ that achieve $\beta_{out}$.  Then $$\beta_{out} = \frac{I(L)+I(R)+|\partial_{out} (L\cup R)|}{|L\cup R|} \geq \frac{\tfrac{1}{d}e(L,L) + \tfrac{1}{d}e(R,R) + \tfrac{1}{d}|\partial(L\cup R)|}{|L\cup R|} = \beta.$$

Now take $L',R'$ that achieve $\beta$.  Then $$\beta = \frac{e(L',L')+e(R',R')+|\partial (L'\cup R')|}{d|L'\cup R'|}\geq \frac{I(L')+I(R')+|\partial_{out}(L'\cup R')|}{d|L'\cup R'|}\geq \tfrac{1}{d}\beta_{out}.$$
\end{proof}

Now, we can prove our main results: first, a bound relating $\beta_{out}$ to $h_{out}$ for all Cayley graphs.  Then we will demonstrate a similar bound relating $\beta$ to $h$.

\begin{theorem} \label{thm_vertexbound}
Let $G$ be a non-bipartite simple Cayley graph corresponding to a group $X$ and a generating set $S$ which satisfies $S^{-1}= S$ and $\text{id}_X\notin S$.  Let $n=|X|$ and $d=|S|$.  Then $$h_{out}(G)\leq 200\beta_{out}.$$
\end{theorem}

\begin{proof}

Choose the disjoint sets $L,R$ that achieve $\beta_{out}$.  Observe that from the definition of $\beta_{out}$, $|\partial_{out} (L\cup R)| \leq \beta_{out}|L\cup R|$.  Set $Y = \{g\in X: \dist_G(g,L\cup R)\geq 2\}$, so that (if $Y$ is non-empty) $\partial_{out} Y = \partial_{out}(L\cup R)$.  Let $\ep>0$ be a constant that we will fix later.  We will first consider the case that $|Y|>\ep n$.

If $\ep n < |Y|<\tfrac{1}{2}n$, then $$h_{out} \leq \frac{|\partial_{out} Y|}{|Y|} \leq \frac{\beta_{out}|L\cup R|}{\ep n}<\frac{1}{\ep}\beta_{out}.$$

If $|Y|\geq \tfrac{1}{2}n$, then $|L\cup R|< \tfrac{1}{2}n$ as $Y$ is disjoint from $L$ and $R$.  And then $$h_{out} \leq \frac{|\partial_{out} (L\cup R)|}{|L\cup R|} \leq \beta_{out}.$$

For the remainder of the proof we will assume $|Y|\leq\ep n$.
For any $g\in X$, define the sets $A(g) = \parens{gL\cap L}\cup \parens{gR\cap R}$ and $B(g) = \parens{gR\cap L}\cup \parens{gL\cap R}$.  Where it is clear we will suppress the input $g$, using $A:=A(g)$ and $B:=B(g)$.  Observe that $A$ and $B$ are disjoint, since $L$ and $R$ are disjoint. We will next bound $|\partial_{out} A|$ and $|\partial_{out} B|$.

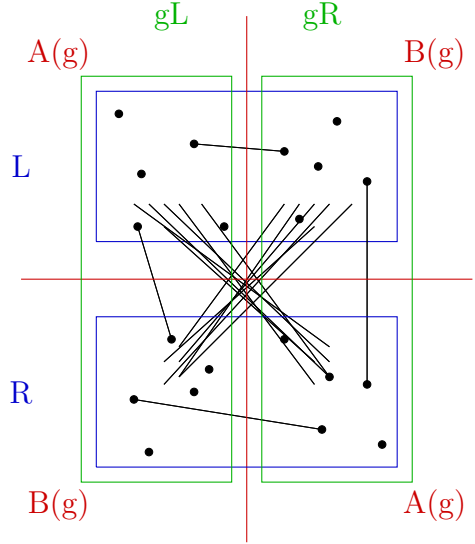
\begin{figure}[!hb]
\begin{center}
\begin{tikzpicture}

\draw[blue] (0,3) -- (4,3)
-- (4,5)-- (0,5)-- (0, 3);

\draw[red] (-.5, 5.5) node{A(g)};
\draw[blue] (-1, 4) node{L};
\draw[green] (1, 6) node{gL};
\draw[red] (4.5, 5.5) node{B(g)};

\draw[black, fill] (3.6, 1.1) circle(0.05cm);
\draw[black, fill] (2.5, 1.7) circle(0.05cm);
\draw[black, fill] (3, .5) circle(0.05cm);
\draw[black, fill] (3.8, .3) circle(0.05cm);
\draw[black, fill] (3.1, 1.2) circle(0.05cm);

\draw[black, fill] (.5,.9) circle(0.05cm);
\draw[black, fill] (.7,.2) circle(0.05cm);
\draw[black, fill] (1,1.7) circle(0.05cm);
\draw[black, fill] (1.3,1) circle(0.05cm);
\draw[black, fill] (1.5,1.3) circle(0.05cm);

\draw[black, fill] (.6, 3.9) circle(0.05cm);
\draw[black, fill] (.55, 3.2) circle(0.05cm);
\draw[black, fill] (.3, 4.7) circle(0.05cm);
\draw[black, fill] (1.7, 3.2) circle(0.05cm);
\draw[black, fill] (1.3, 4.3) circle(0.05cm);

\draw[black, fill] (2.5, 4.2) circle(0.05cm);
\draw[black, fill] (2.7, 3.3) circle(0.05cm);
\draw[black, fill] (3.6, 3.8) circle(0.05cm);
\draw[black, fill] (2.95, 4) circle(0.05cm);
\draw[black, fill] (3.2, 4.6) circle(0.05cm);

\draw[black, fill] (3.1, 3.5) -- (.9, 1.1);
\draw[black, fill] (2.9, 3.5) -- (1.1, 1.4);
\draw[black, fill] (2.7, 3.5) -- (1.1, 1.2);
\draw[black, fill] (2.5, 3.5) -- (1.1, 1.6);
\draw[black, fill] (2.9, 3.2) -- (0.9, 1.4);
\draw[black, fill] (3.4, 3.5) -- (1.1, 1.2);

\draw[black, fill] (1.1, 3.5) -- (2.9, 1.1);
\draw[black, fill] (0.9, 3.5) -- (3.1, 1.4);
\draw[black, fill] (0.7, 3.5) -- (3.1, 1.2);
\draw[black, fill] (0.5, 3.5) -- (3.1, 1.6);
\draw[black, fill] (0.9, 3.2) -- (2.9, 1.4);
\draw[black, fill] (1.4, 3.5) -- (3.1, 1.2);

\draw[black, fill] (2.5, 4.2) -- (1.3, 4.3);
\draw[black, fill] (.5,.9)  -- (3, .5);
\draw[black, fill] (3.6, 3.8) -- (3.6, 1.1) ;
\draw[black, fill] (.55, 3.2) -- (1,1.7) ; 

\draw[red] (-1, 2.5) -- (5, 2.5);
\draw[red] (2, 6) -- (2, -1);

\draw[blue] (0,0) -- (4,0)
-- (4,2)-- (0,2)-- (0,0);

\draw[red] (-.5, -.5) node{B(g)};
\draw[blue] (-1, 1) node{R};
\draw[green] (3, 6) node{gR};
\draw[red] (4.5, -.5) node{A(g)};

\draw[green] (2.2,-0.2) rectangle (4.2, 5.2);
\draw[green] (-0.2,-0.2) rectangle (1.8, 5.2);
\end{tikzpicture}

\end{center}
\caption{Illustration of $gL$, $gR$, and $A(g)$, $B(g)$.}
\end{figure}

Consider the set $(\partial_{out}A)\cap B$: any vertex in that set must be counted by one of $I(L), I(R), I(gL)$, or $I(gR)$.  Any other vertex in $\partial_{out}A$ must be in $\partial_{out}(L\cup R)$ or $\partial_{out}((L\cup R)g)$.  By symmetry the same holds for $\partial_{out}B$.
It follows that $|\partial_{out}A(g)|,|\partial_{out}B(g)|\leq 2\parens{I(L)+I(R)+|\partial_{out}(L\cup R)|}$, and therefore $$h_{out}(G)\leq \frac{2\parens{I(L)+I(R)+|\partial_{out}(L\cup R)|}}{\min |A(g)|,|B(g)|}.$$  Here we use the facts that $I(gL) = I(L)$ and $I(gR) = I(R)$.  It is simple to see that the numerator is $2\beta_{out}|L\cup R|$, it remains to bound $\min\{|A(g)|,|B(g)|\}$.

For this step, we will use a technique developed in \cite{Freiman}, which was used to prove similar results in \cite{Biswas,BS,BGGT}.

Without loss of generality, assume that $|L| \geq |R|$.  Notice that $|L\cup R| = n-|Y|-|\partial_{out} (L\cup R)|\geq n-\ep n -\beta_{out}n$, and hence by assumption $|L|\geq \tfrac{1-\ep -\beta_{out}}{2}n$.  

Suppose that $|L|\geq \tfrac{1+\ep}{2}n$.  Observe that $|L|-I(L)\leq |L^c|$ and so $I(L) \geq \ep n$.  Also observe that as a general bound $h_{out}\leq \tfrac{n+1}{n-1}\leq 2$.  In this case we have the bound $\beta_{out}\geq \frac{I(L)}{n} \geq \ep$, and so $h_{out}\leq 2 \leq \tfrac{2}{\ep}\beta_{out}$.  We will consider the other case, where $\tfrac{1-\ep -\beta_{out}}{2}n\leq |L|\leq \tfrac{1+\ep}{2}n$.

Assume for contradiction that there is no element $g\in X$ for which $|L\cap gL|\in (\delta |L|,(1-\delta)|L|)$, where $\delta$ is a constant we will define later.  Define the sets $X_1 = \{g:|L\cap gL|\geq (1-\delta)|L|\}$ and $X_2 = \{g:|L\cap gL|\leq \delta |L|\}$.  By assumption $X_1,X_2$ is a partition of $X$.

We can show that if $\delta < \tfrac{1}{3}$, then $X_1^2 \subset X_1$.  Let $g,h\in X_1$, then 
\begin{eqnarray*}
|ghL\cap L| & \geq &|L|-|gL-ghL|-|L-gL|\\
&\geq &(1-2\delta)|L| > \delta |L|\,,
\end{eqnarray*} 
 and so assuming $\delta < \tfrac{1}{3}$, $gh\in X_1$.

Similarly we can show that if $\delta < 1-\tfrac{2}{3(1-\ep-\beta_{out})}$, then $X_2^2 \subset X_1$.  Let $g,h\in X_2$, then \begin{eqnarray*}
|ghL \cap L| & \geq &|L-gL| - |X-(gL\cap ghL)|\\
& = & 3|L|-n-|L\cap gL|-|gL\cap ghL|\\ 
& \geq &(3-\tfrac{2}{1-\ep-\beta_{out}})|L|- 2\delta |L| > \delta|L|\,,\end{eqnarray*}
 and so assuming $\delta < 1-\tfrac{2}{3(1-\ep-\beta_{out})}$, $gh\in X_1$.

Because $X_1^2\subset X_1$, $X_1$ is a subgroup of $X$.  Suppose that $X_2$ is empty, so that $X_1 = X$.  Then if $\delta<\tfrac{1}{3}$, $|L|^2 = \sum_g |gL\cap L| > \tfrac{2}{3}|L|n$; i.e., $|L|> \tfrac{2}{3}n$.  We are already assuming that $|L|\leq \tfrac{1+\ep}{2}n$, so there is a contradiction as long as we eventually choose $\ep < \tfrac{1}{3}$.  

On the contrary, we have that $X_2$ is non-empty.  As $X_1$ is a proper subgroup of $X$ and $X_2$ is its complement with $X_2^2\subset X_1$, it follows that $X_1$ is a subgroup of index $2$ with $X_2$ as its unique non-trivial coset.
Because of this,
\begin{eqnarray*}
\tfrac{n}{2}(1-\delta)|L| & \leq &\sum_{g\in X_1} |L\cap gL| = |L\cap X_1|^2 + |L\cap X_2|^2\\ & = &|L|^2 - 2|L\cap X_1||L\cap X_2|\\ 
&\leq  &\tfrac{n}{2}(1+\ep)|L|-2|L\cap X_1||L\cap X_2|\,.
\end{eqnarray*}

It follows that $$|L\cap X_1||L\cap X_2|\leq (\ep + \delta)\tfrac{n}{4}|L|.$$

This means that there is $i \in \{1,2\}$ so that $|L\cap X_i|\leq \sqrt{(\ep + \delta)\tfrac{n}{4}|L|}$.  Let $X_j$ be the other coset of $X_1$, for which $|L\cap X_j| = |L|-|L\cap X_i| \geq |L| - \sqrt{(\ep + \delta)\tfrac{n}{4}|L|}$.

Let $s\in S$, consider the set $X_j s \cap X_j$.  If $g\in X_j s \cap X_j$, either (1) $g\in X_j-L$, (2) $g\in (X_j-L)s$ or (3) $\{g,gs^{-1}\}$ is an edge in $E(L,L)$.  \begin{eqnarray*}
|X_j s\cap X_j| &\leq &2|X_j-L|+I(L) \leq  2\parens{\tfrac{n}{2}-|X_j \cap L|}+\beta_{out}n \\  
&\leq &2\parens{\tfrac{n}{2}-|L|+\sqrt{(\ep + \delta)\tfrac{n}{4}|L|}}+\beta_{out}n \\
&\leq &2\parens{\tfrac{n}{2}(\ep + \beta_{out}) + \tfrac{n}{2}\sqrt{\tfrac{1}{2}(\ep + \delta)(1+\ep)}}+\beta_{out} n\\
& = & \parens{\ep + 2\beta_{out} + \sqrt{\tfrac{1}{2}(\ep + \delta)(1+\ep)}}n.
\end{eqnarray*}

Assuming we choose $\delta < (\tfrac{1}{2}-\ep -2\beta_{out})^2\tfrac{2}{1+\ep}-\ep$, then $|X_js \cap X_j|< \tfrac{n}{2}$.  But, if $s\in X_1$, then $|X_js \cap X_j| = |X_j| = \tfrac{n}{2}$.  Therefore $S\subset X_2$.  If $g\in X_1$ and $s\in S$, then $gs\in gX_2 = X_2$.  Likewise if $g\in X_2$ and $s\in S$, then $gs\in gX_2 = X_1$.  Every edge of $G$ is incident to one vertex from $X_1$ and one in $X_2$, in other words, $G$ is bipartite.  This is our desired contradiction.

So instead, let $g$ be an element of $X$ for which $|gL\cap L|\in (\delta |L|,(1-\delta)|L|)$.  $$|A(g)| \geq |gL \cap L| \geq \delta |L|\geq \frac{\delta (1-\ep -\beta_{out})n}{2}.$$  $$|B(g)|\geq |gL \cap R|\geq |L-gL|-|(L\cup R)^c|\geq \frac{\delta(1-\ep-\beta_{out})n}{2} - (\ep + \beta_{out})n.$$

Now we can complete the bound on the the vertex expansion: $$h_{out}(G)\leq \frac{2\beta_{out} n}{\min |A(g)|,|B(g)|}\leq \frac{2\beta_{out} n}{\frac{\delta(1-\ep-\beta_{out})n}{2} - (\ep + \beta_{out})n} = \frac{2\beta_{out}}{\tfrac{\delta(1-\ep-\beta_{out})}{2} - (\ep + \beta_{out})}.$$

At this point we will make our choices of $\delta$ and $\ep$.  We have previously required that $$\delta < \tfrac{1}{3},\text{ } \delta < 1-\tfrac{2}{3(1-\ep-\beta_{out})}, \text{ and } \delta < (\tfrac{1}{2}-\ep -2\beta_{out})^2\tfrac{2}{1+\ep}-\ep.$$

If $\beta_{out}\geq \ep$, then we have a bound $h_{out} \leq 2 \leq \tfrac{2}{\ep}\beta_{out}$.  In the other case, we assume that $\beta_{out}\leq \ep$; it is now enough to require that $$\delta < \tfrac{1}{3},\text{ } \delta < 1-\tfrac{2}{3-6\ep},\text{ and } \delta < (\tfrac{1}{2}-3\ep)^2\tfrac{2}{1+\ep}-\ep.$$  To satisfy these restrictions we will set $\ep = \tfrac{1}{100}$ and $\delta = \tfrac{1}{5}$.
To summarize all the cases,
$$h_{out}\leq \begin{cases}
\tfrac{1}{\ep}\beta_{out} & \text{if } \ep n \leq |Y| \leq \tfrac{1}{2}n,\\
\beta_{out} & \text{if } |Y| \geq \tfrac{1}{2}n,\\
\tfrac{2}{\ep}\beta_{out} & \text{if }\max |L|,|R|\geq (\tfrac{1+\ep}{2})n,\\
\tfrac{2}{\ep}\beta_{out} & \text{if }\beta_{out}\geq \ep,\\
\frac{2\beta_{out}}{\tfrac{\delta(1-2\ep)}{2} - 2\ep} & \text{otherwise.}
\end{cases}$$

Substituting our choices of $\ep$ and $\delta$ into each of these bounds completes the proof.

\end{proof}

In this theorem, we derive a similar relationship between the edge versions of $h$ and $\beta$.  The proof of this result is almost identical to that of Theorem \ref{thm_vertexbound}.

\begin{theorem} \label{thm_edgebound}
Let $G$ be a non-bipartite simple Cayley graph corresponding to a group $X$ and a generating set $S$ which satisfies $S^{-1}= S$ and $\text{id}_X\notin S$.  Let $n=|X|$ and $d=|S|$.  Assume that $\beta(G) < \tfrac{1}{15d}$.  Then $$h\leq \frac{100 \beta}{1-15d\beta}.$$
\end{theorem}

\begin{proof}

Choose the disjoint sets $L,R$ that achieve $\beta$, so that $|\partial (L\cup R)| \leq d\beta|L\cup R|$.  Set $Y = X-(L\cup R)$, so that $\partial Y = \partial(L\cup R)$.  Let $\ep>0$ be a constant that we will fix later.  We will first consider the case that $|Y|>\ep n$.

If $\ep n < |Y|<\tfrac{1}{2}n$, then $$h \leq \frac{|\partial Y|}{d|Y|} \leq \frac{d\beta|L\cup R|}{d\ep n}<\frac{1}{\ep}\beta.$$

If $|Y|\geq \tfrac{1}{2}n$, then $|L\cup R|< \tfrac{1}{2}n$ as $Y$ is disjoint from $L$ and $R$.  And then $$h \leq \frac{|\partial (L\cup R)|}{d|L\cup R|} \leq \beta.$$

For the remainder of the proof we will assume $|Y|\leq\ep n$.
For any $g\in X$, define the sets $A(g) = \parens{gL\cap L}\cup \parens{gR\cap R}$ and $B(g) = \parens{gR\cap L}\cup \parens{gL\cap R}$.

Observe that $A$ and $B$ are disjoint, since $L$ and $R$ are disjoint.  Also define $Z(g) = Y\cup gY$ to be the complement of $A(g)\cup B(g)$, so that $\{A,B,Z\}$ is a partition of $X$.  Note that $\partial A = E(A,B)\cup E(A,Z)$, similarly $\partial B = E(B,A)\cup E(B,Z)$, and 
so we next bound $|\partial_{out} A|$ and $|\partial_{out} B|$.  

\medskip
Consider the set $E(A,B)$: if $(v,w)\in A\times B$, then $v$ and $w$ must be common members of one of $L,R,gL,gR$.  If $(v,w)\in E(A,B)$, then $(v,w)$ must be in one of the following four sets: $E(L,L),E(R,R),E(gL,gL),E(gR,gR)$.  And so $e(A,B)\leq e(L,L)+e(R,R)+e(gL,gL)+e(gR,gR) = 2\parens{e(L,L)+e(R,R)}$.  

\smallskip
Similarly, any edge in $E(A,Z)$ or $E(B,Z)$ must be in $E(L\cup R,Y)$ or $E (g(L\cup R),gY)$, and so $e(A,Z),e(B,Z)\leq |\partial (L\cup R)|+|\partial g(L\cup R)| = 2|\partial (L\cup R)|$.

\medskip
It follows that $|\partial A(g)|,|\partial B(g)|\leq 2\parens{e(L,L)+e(R,R)+|\partial (L\cup R)|}$, and therefore $$h(G)\leq \frac{2\parens{e(L,L)+e(R,R)+|\partial (L\cup R)|}}{d\cdot\min\{ |A(g)|,|B(g)|\}}.$$  Here we use the facts that $e(gL,gL) = e(L,L)$ and $e(gR,gR) = e(R,R)$.  By our choice of $L$ and $R$, the numerator is precisely $2d\beta |L\cup R|$, hence it remains to bound $\min \{|A(g)|,|B(g)|\}$.

\medskip
Once again we will use a technique of \cite{Freiman}.

Without loss of generality, assume that $|L| \geq |R|$.  Notice that $|L\cup R| = n-|Y|\geq n-\ep n$, and hence by assumption $|L|\geq \tfrac{1-\ep}{2}n$.

Suppose that $|L|\geq \tfrac{1+\ep}{2}n$.  Observe that $e(L,R^C)\geq d(|L|-|R|)\geq d\ep n$.  In this case we have the bound $\beta\geq \frac{e(L,R^C)}{nd} \geq \ep$, and so $h\leq 1 \leq \tfrac{1}{\ep}\beta$.

\smallskip
We will consider the other case, where $\tfrac{1-\ep}{2}n\leq |L|\leq \tfrac{1+\ep}{2}n$.

Assume for contradiction that there is no element $g\in X$ for which $|L\cap gL|\in (\delta |L|,(1-\delta)|L|)$, where $\delta$ is a constant we will define later.  Define the sets $X_1 = \{g:|L\cap gL|\geq (1-\delta)|L|\}$ and $X_2 = \{g:|L\cap gL|\leq \delta |L|\}$.  By assumption $X_1,X_2$ is a partition of $X$.

We can show that if $\delta < \tfrac{1}{3}$, then $X_1^2 \subset X_1$.  Let $g,h\in X_1$, then
\begin{eqnarray*}
|ghL\cap L| & \geq &|L|-|gL-ghL|-|L-gL|\\ 
&\geq &(1-2\delta)|L| > \delta |L|\,,\end{eqnarray*} 
and so assuming $\delta < \tfrac{1}{3}$, $gh\in X_1$.

Similarly we can show that if $\delta < 1-\tfrac{2}{3-3\ep}$, then $X_2^2 \subset X_1$.  Let $g,h\in X_2$, then 
\begin{eqnarray*}
|ghL \cap L| & \geq &|L-gL| - |X-(gL\cap ghL)|\\ & = & 3|L|-n-|L\cap gL|-|gL\cap ghL|\\ 
& \geq &(3-\tfrac{2}{1-\ep})|L| - 2\delta |L|  > \delta|L|\,,
\end{eqnarray*}
and so assuming $\delta < 1-\tfrac{2}{3-3\ep}$, $gh\in X_1$.

Because $X_1^2\subset X_1$, $X_1$ is a subgroup of $X$.  Suppose that $X_2$ is empty, so that $X_1 = X$.  Then if $\delta <\tfrac{1}{3}$, $|L|^2 = \sum_g |gL\cap L| > \tfrac{2}{3}|L|n$; i.e., $|L|> \tfrac{2}{3}n$.  We already assumed that $|L|\leq \tfrac{1+\ep}{2}n$, so there is a contradiction as long as we eventually choose $\ep < \tfrac{1}{3}$.  

On the contrary, we have that $X_2$ is non-empty.  As $X_1$ is a proper subgroup of $X$ and $X_2$ is its complement with $X_2^2\subset X_1$, it follows that $X_1$ is a subgroup of index $2$ with $X_2$ as its unique non-trivial coset.  As a result,
\begin{eqnarray*}
\tfrac{n}{2}(1-\delta)|L| & \leq &\sum_{g\in X_1} |L\cap gL|\\ 
& = &|L\cap X_1|^2 + |L\cap X_2|^2\\ & = &|L|^2 - 2|L\cap X_1||L\cap X_2|\\ 
& \leq  &\tfrac{n}{2}(1+\ep)|L|-2|L\cap X_1||L\cap X_2|\,.\end{eqnarray*}

It follows that $$|L\cap X_1||L\cap X_2|\leq (\ep + \delta)\tfrac{n}{4}|L|.$$

This means that $\exists i \in \{1,2\}$ so that $|L\cap X_i|\leq \sqrt{(\ep + \delta)\tfrac{n}{4}|L|}$.  Let $X_j$ be the other coset of $X_1$, for which $|L\cap X_j| = |L|-|L\cap X_i|\geq |L| - \sqrt{(\ep + \delta)\tfrac{n}{4}|L|}$.

Let $s\in S$, consider the set $X_j s \cap X_j$.  If $g\in X_j s \cap X_j$, then either (1) $g\in X_j-L$, or (2) $g\in (X_j-L)s$ or (3) $g\in I(L)$.  \begin{eqnarray*}
|X_j s\cap X_j| & \leq &2|X_j-L|+I(L) \\ 
&\leq  & 2\parens{\tfrac{n}{2}-|X_j \cap L|}+d\beta n \\ 
& \leq & 2\parens{\tfrac{n}{2}-|L|+\sqrt{(\ep + \delta)\tfrac{n}{4}|L|}}+d\beta n \\
&\leq &2\parens{\tfrac{\ep n}{2} + \tfrac{n}{2}\sqrt{
\tfrac{1}{2}(\ep + \delta)(1+\ep)}}+d\beta n\\
& = & \parens{\ep + d\beta + \sqrt{\tfrac{1}{2}(\ep + \delta)(1+\ep)}}n\,.
\end{eqnarray*}

Assuming that we chose $\delta < (\tfrac{1}{2}-\ep -d\beta )^2\tfrac{2}{1+\ep}-\ep$, then $|X_js \cap X_j|< \tfrac{n}{2}$.  But, if $s\in X_1$, then $|X_js \cap X_j| = |X_j| = \tfrac{n}{2}$.  Therefore $S\subset X_2$.  If $g\in X_1$ and $s\in S$, then $gs\in gX_2 = X_2$.  Likewise if $g\in X_2$ and $s\in S$, then $gs\in gX_2 = X_1$.  Every edge of $G$ is incident to one vertex from $X_1$ and one in $X_2$, in other words, $G$ is bipartite.  This is our desired contradiction.

So instead, let $g$ be an element of $X$ for which $|gL\cap L|\in (\delta |L|,(1-\delta)|L|)$.  $$|A(g)| \geq |gL \cap L| \geq \delta |L|\geq \frac{\delta (1-\ep )n}{2}.$$  $$|B(g)|\geq |gL \cap R|\geq |L-gL|-|Y|\geq \frac{\delta(1-\ep)n}{2} - \ep n.$$

Now we can complete the bound of the Cheeger constant; $$h(G)\leq \frac{2d\beta n}{d\cdot\min |A(g)|,|B(g)|}\leq \frac{2d\beta n}{\frac{\delta(1-\ep)n}{2} - \ep n} = \frac{2\beta}{\tfrac{\delta(1-\ep)}{2} - \ep}.$$

At this point we will make our choices of $\delta$ and $\ep$.  We have previously required that $$\delta < \tfrac{1}{3},\text{ } \delta < 1-\tfrac{2}{3-3\ep},\text{ and } \delta < (\tfrac{1}{2}-\ep -2d\beta)^2\tfrac{2}{1+\ep}-\ep.$$

To satisfy these bounds, take $\ep = \tfrac{1}{100}$ and $\delta = \tfrac{32}{100}-4d\beta$.

$$h(G)\leq \frac{2\beta}{\parens{\tfrac{32}{100}-4d\beta}\tfrac{99}{200}-\tfrac{1}{100}}< \frac{2\beta}{\tfrac{14}{100}-2d\beta}<\frac{15\beta}{1-15d\beta}.$$

After using our values of $\ep$ and $\delta$, we can summarize all the cases:
$$h\leq \begin{cases}
100 \beta & \text{if } \tfrac{1}{100} n \leq |Y| \leq \tfrac{1}{2}n,\\
\beta & \text{if } |Y| \geq \tfrac{1}{2}n,\\
100 \beta& \text{if }\max |L|,|R|\geq (\tfrac{101}{200})n,\\
\frac{15\beta}{1-15d\beta} & \text{otherwise.}
\end{cases}$$

As a general bound, we may conclude with $$h\leq \frac{100 \beta}{1-15d\beta}\,,$$  completing the proof.

\end{proof}

Now by combining our results with Trevisan's bound, we obtain the following theorem that improves on the main result of \cite{Biswas}.

\begin{theorem}\label{thm_finalbound}
Let $G$ be a $d$-regular Cayley graph.  
\begin{enumerate}
\item There is a universal constant $C_1$ so that $$\mu \geq \frac{C_1h_{out}^2}{d^2}.$$
\item There is a universal constant $C_2$ so that the following holds: if $\beta \leq \frac{1}{30d}$ then $$\mu \geq C_2 h^2.$$  Otherwise, $$\mu \geq  \frac{C_2}{d^2}.$$\end{enumerate}
\end{theorem}

\begin{proof}
To see Item 1, recall from Theorem \ref{thm_trevisan} that $\mu \geq \tfrac{1}{2}\beta^2$.  From Theorem \ref{thm_betabounds} recall that $\beta\geq \tfrac{1}{d}\beta_{out}$.  And our result in Theorem \ref{thm_vertexbound} we have $\beta_{out}\geq \tfrac{1}{200}h_{out}$.  Combining these inequalities gives the desired result, with $C_1 = \tfrac{1}{80,000}$.

To see Item 2, we again use $\mu\geq \tfrac{1}{2}\beta^2$.  The result of Theorem \ref{thm_edgebound} is that $h \leq \frac{100\beta}{1-15d\beta}$.  If $\beta \leq \frac{1}{30d}$, then $h \leq 200\beta$.  On the other hand, if $\beta > \frac{1}{30d}$, then we see directly that $\mu \geq \frac{1}{2\cdot(30d)^2}$.  These two cases give the desired result, with $C_2 = \tfrac{1}{80,000}$.

\end{proof}

\subsection{Examples}

As an illustrative example of our proof method in Theorems \ref{thm_vertexbound} and \ref{thm_edgebound}, we will examine an odd cycle.

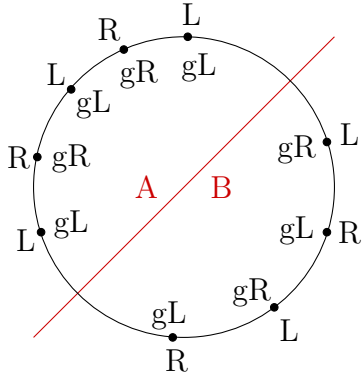
\begin{figure}[!hb]
\begin{center}
\begin{tikzpicture}

\draw[red] (-2,-2) -- (2,2);
\draw (0,0) circle (2cm);

\draw (.1, 2.3) node{L}; 
\draw (.2,1.6) node{gL};
\draw[black, fill] (.05,2) circle(0.05cm);

\draw (-1, 2.1) node{R}; 
\draw (-.6,1.5) node{gR};
\draw[black, fill] (-0.8,1.83) circle(0.05cm);

\draw (-1.7, 1.5) node{L}; 
\draw (-1.2,1.1) node{gL};
\draw[black, fill] (-1.5, 1.3) circle(0.05cm);

\draw (-2.2, .4) node{R}; 
\draw (-1.5, .38) node{gR};
\draw[black, fill] (-1.95,.4) circle(0.05cm);

\draw (-2.1, -.7) node{L}; 
\draw (-1.5, -.5) node{gL};
\draw[black, fill] (-1.9,-0.6) circle(0.05cm);

\draw[red](-.5,0) node{A};
\draw[red](.5,0) node{B};

\draw (2.2, .7) node{L}; 
\draw (1.5, .5) node{gR};
\draw[black, fill] (1.9,0.6) circle(0.05cm);

\draw (2.2, -.6) node{R}; 
\draw (1.5, -.58) node{gL};
\draw[black, fill] (1.9,-0.6) circle(0.05cm);

\draw (1.4, -1.9) node{L}; 
\draw (.9, -1.4) node{gR};
\draw[black, fill] (1.2,-1.6) circle(0.05cm);

\draw (-.1, -2.3) node{R}; 
\draw (-.2, -1.7) node{gL};
\draw[black, fill] (-.15, -2) circle(0.05cm);

\end{tikzpicture}
\caption{Illustration of $A$ and $B$ in the odd cycle $C_9$.}
\end{center}
\end{figure}

\begin{example}\label{example_cycle}
Let $X = Z_{2k+1}$ where $k$ is a positive integer, and consider the Cayley graph $G = (X,\{\pm 1\})$.  Using the methods of our theorem, we can use the $\beta_{out}(G)$-achieving almost-bipartition to obtain a vertex cut which approximates $h_{out}(G)$.  Similarly we can use the $\beta(G)$-achieving almost-bipartition to obtain a cut which approximates $h(G)$.
\end{example}

\begin{proof}

First we solve the vertex-expansion problem.
$\beta_{out} = \tfrac{1}{2k}$, achieved by taking $L$ to be the odd integers $\{1,3,\dots, 2k-1\}$, $R = \{2,4,\dots, 2k\}$ and leaving $2k+1$ uncolored.

Take $g = k$, so that $gL = \{k+1,k+3,\dots, k-2\}$ and $gR = \{k+2,k+4,\dots, k-1\}$.  $Z(k) = \{0,k\}$.  If $k$ is even, then $L\cap gL = \{k+1,k+3,\dots, 2k-1\}$ and $R\cap gR = \{k+2,k+4,\dots ,2k\}$, thus $A(k) = \{k+1,\dots,2k\}$ and $B(k) = \{1,\dots,k-1\}$.  On the other hand if $k$ is odd, it can be seen that $A(k) = \{1,\dots,k-1\}$ and $B(k) = \{k+1,\dots,2k\}$.  In either case, $\partial_{out} A(k) = \partial_{out} B(k) = Z(k)$ and we have the bound $$h_{out}\leq \frac{|Z(k)|}{\min |A(k)|,|B(k)|} = \frac{2}{k-1}.$$  Observe that in this case we have approximately achieved the actual value $h_{out} = \tfrac{2}{k}$ by choosing an optimal value of $k$.  It is well known that $\mu = \Theta(\frac{1}{k^2})$, and clearly $d = \Theta(1)$, so our bound $\mu \gtrsim \frac{h_{out}^2}{d^2}$ is tight up to a constant factor.

Now, we use similar methods to work with the edge expansion.
$\beta = \tfrac{1}{2k+1}$, achieved by taking $L$ to be the odd integers $\{1,3,\dots, 2k+1\}$, $R = \{2,4,\dots, 2k\}$.

Take $g = k$, so that $gL = \{k+1,k+3,\dots, k-2,k\}$ and $gR = \{k+2,k+4,\dots, k-1\}$.  $Z(k) = \{\}$  If $k$ is even, then $L\cap gL = \{k+1,k_3,\dots,2k+1\}$ and $R\cap gR = \{k+2,k+4,\dots,2k\}$, thus $A(k) = \{k+1,\dots,2k+1\}$ and $B(k) = \{1,\dots, k\}$.  On the other hand if $k$ is odd, it can be seen that $A(k) = \{1,\dots, k\}$ and $B(k) = \{k+1,\dots,2k+1\}$.  In either case, $\partial A(k) = \partial B(k)= \{\{k,k+1\},\{2k+1,1\}\}$ and we have the bound $$h\leq \frac{2}{2\min |A(k)|,|B(k)|} = \frac{1}{k-1}.$$  Observe that in this case we have approximately achieved the actual value $h = \tfrac{1}{k}$ by choosing an optimal value of $k$.  It is well known that $\mu = \Theta(\frac{1}{k^2})$, so our bound $\mu \gtrsim h^2$ is tight up to a constant factor.

\end{proof}

\begin{remark}
Notice that in this example, we did not need to use Freiman's method, as it is simple to explicitly find a value $g$ for which $A(g)$ and $B(g)$ are both $\Theta(n)$.
\end{remark}

We will now give a simple example that shows our bound $\beta \gtrsim h$ is not tight for general Cayley graphs, and indeed that there cannot be a reverse inequality of the form $h\gtrsim f(\beta,d)$ for any non-trivial function $f$.

\begin{example}\label{example_no_converse}
Let $X = Z_3\times Z_{2k+1}$ where $k$ is a positive integer, let $S = \{(\pm 1,0),(0,\pm 1)\}$, and let $G = (X,S)$ be the Cayley graph.  Then $h_{out}(G) \ll \beta_{out}(G).$
\end{example}

\begin{proof}
Consider the set $A = \{[k]\times Z_3\}$, with $|\partial_{out}(A)| = 6$ and $|A| = 3k$.  Using $A$ as a candidate we see that $h_{out} \leq \frac{2}{k}$.

Let $L,R$ be a candidate bipartition of $X$.  For any $3$-cycle $C$ in $G$, if $L$ or $R$ intersects $C$ then at least one vertex of $C$ must be in $\partial_{out}(L\cap R)$ or be counted by $I(L)$ or $I(R)$.

That means that \begin{eqnarray*}b_{out}(L,R) & = &\frac{I(L)+I(R)+|\partial_{out}(L\cup R)}{|L\cup R|}\\
& \geq &\frac{\sum_C I_C(L)+I_C(R)+|\partial_{out, C}(L\cup R)|}{|L\cup R|}\\ 
&\geq &\frac{\tfrac{1}{3}\sum_C|(L\cup R)\cap C|}{|L\cup R|} = \frac{1}{3}\,.\end{eqnarray*}

So $\beta_{out} \geq \tfrac{1}{3} \gg \frac{2}{k}\geq h_{out}$.

\end{proof}

At a high level we are looking for bounds on $h$ in terms of $\beta$ and $d$.  As $\beta$ and $d$ are both $\Theta(1)$, this example tells us that there can be no lower bound on $h$ that applies to all Cayley graphs.  Observe that a similar analysis gives a similar result for $h$ and $\beta$ on the same graph.

\section{Open Questions}

\begin{itemize}
    \item 
    Recall the Cheeger inequalities $2h\geq \lambda \geq \tfrac{1}{2}h^2$, where $\lambda:=1-\lambda_2$.  A problem of general interest is to categorize the graphs for which $\lambda\approx h$ and those for which $\lambda \approx h^2$.  Similarly, we can ask if there are some non-trivial classes of non-bipartite graphs for which $\mu \approx \beta^2$ (or alternately $\mu \approx \beta$).
    In particular, there has recently been investigation into various definitions of the discrete curvature.  For example, Klartag et al.~\cite{KKRT} (see also \cite{EF}) demonstrated that if a graph has non-negative curvature in the sense of the curvature-dimension inequality, then $16dh^2 \geq \lambda$; that is, $\lambda \approx h^2$.  A class of graphs for which this curvature bound holds is Cayley graphs of abelian groups.  Is there a definition of discrete curvature that permits the characterization of a class of graphs for which $\mu \approx \beta^2$?
    \item 
    In our result $\mu \geq \frac{Ch_{out}}{d^2}$ our focus was on obtaining the correct dependence of $\mu$ on $h_{out}$ and we did not explore the tightness in terms of degree $d$.  In the proof we first relate $\mu$ to $\beta$ and then use the simple bound $\beta \geq \tfrac{1}{d}\beta_{out}$ from Theorem~\ref{thm_betabounds}.  Bobkov, Houdr\'{e}, and the third author~\cite{BHT}  introduced a functional constant $\lambda_{\infty}$ and used the proof methods of Cheeger inequalities to demonstrate an analogous relationship between $\lambda_{\infty}$ and $h_{out}$\,.
   
    Is it possible to do the same for $\beta_{out}$; that is, can we define a functional constant $\mu_{\infty}$ (say) and prove directly a relationship between $\mu_{\infty}$ and $\beta_{out}$.  This would be in contrast to our current proof which relates $\mu$ and $\beta_{out}$, using $\beta$ as an intermediary.
    
    \item 
    Biswas and Saha~\cite{BS} proved that for any non-bipartite Cayley sum graph (that is, a graph defined by the relation $(g,h)\in E$ iff $gh\in S$ for some generating set $S$), $\mu \geq \frac{Ch_{out}^4}{d^8}$ for a universal constant $C$.  To obtain this result they modified the proof method of Biswas's similar result for Cayley graphs~\cite{Biswas}.  The modification is necessary because the original result makes use of the vertex-transitivity of a Cayley graph; a Cayley sum graph need not be transitive.  Is it possible to extend our Theorem~\ref{thm_finalbound} to the setting of Cayley sum graphs in a similar way?
\end{itemize}

%
%
%
%
%
\bibliography{biblio}
\bibliographystyle{plain}

\end{document}